# On the propagation of singularities of constant curvature, convex hypersurfaces

19 January 2024

Graham Smith[*]

*Dedicated to Marcos Dajczer on the occasion of his 75th birthday*

**Abstract:** We describe the structure of the singular sets of constant curvature, convex hypersurfaces in hyperbolic space for general convex curvature functions. We apply this result to the study of the ideal Plateau problem in hyperbolic space for such curvature functions.

**1 - Introduction.** In the study of certain non-linear partial differential equations, we are in the fortunate position of having complete geometric descriptions of the singular sets of their solutions, which in certain favourable situations even allows us to eliminate these singularities entirely. This is the case, for example, for equations of real Monge-Ampère type (see, for example, [2], [11], and [18]), where the singular sets are unions of convex hulls of closed subsets of the boundary of the domain of the function in question, as in Figure 1.1. In this note, we study the singular sets of convex hypersurfaces of constant curvature for a general class of curvature functions to be specified presently. The corresponding partial differential equations are of hessian type, and we will see that the general structure of the singular sets of their solutions is similar to that just described. We will work in hyperbolic space, although analogous results can also be shown to hold in the sphere and in euclidean space.

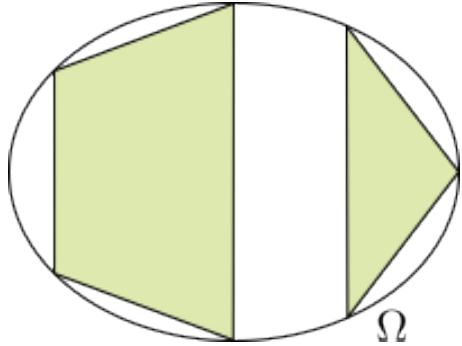

**Figure 1.1 - The structure of the singular set** - Given a solution of the real Monge-Ampère equation defined in some convex domain $\Omega$, its singular set is a union of convex hulls of closed subsets of the boundary of this domain.

In order to state our main result, a few definitions will be required. Let $\mathbb{H}^{m+1}$ denote $(m+1)$-dimensional hyperbolic space, let $\partial_\infty\mathbb{H}^{m+1}$ denote its ideal boundary. Let Conv denote the space of closed, convex subsets of $\mathbb{H}^{m+1} \cup \partial_\infty\mathbb{H}^{m+1}$ furnished with the Hausdorff topology, and let $\text{Conv}^o$ denote the subspace consisting of those closed, convex subsets which have non-trivial interior.

Let $K$ be a convex curvature function, as defined in Section 3. Given an open subset $\Omega$ of $\mathbb{H}^{m+1}$, and a positive real number $\kappa$, let $\text{Conv}^o_\kappa(\Omega)$ denote the space of closed, convex subsets of $\mathbb{H}^{m+1} \cup \partial_\infty\mathbb{H}^{m+1}$, with non-trivial interior, such that $\partial X \cap \Omega$ is smooth with constant $K$-curvature equal to $\kappa$ (see Figure 1.2).

Let $\partial\text{Conv}^o_\kappa(\Omega)$ denote the topological boundary of $\text{Conv}^o_\kappa(\Omega)$ in $\text{Conv}^o$. Given $X \in \partial\text{Conv}^o_\kappa(\Omega)$ and $x \in \partial X \cap \Omega$, we say that $X$ is *regular* at $x$ whenever $X \in \text{Conv}^o_\kappa(\Omega')$ for some smaller neighbourhood $\Omega'$ of $x$, and we say that $X$ is *singular* at this point otherwise. Let $\text{Reg}(X,\Omega)$ denote the set of points of $\partial X \cap \Omega$ at which $X$ is regular, and let $\text{Sing}(X,\Omega)$ denote the set of points at which it is singular.

---

[*] Departamento de Matemática, Pontifícia Universidade Católica, Rio de Janeiro, Brazil

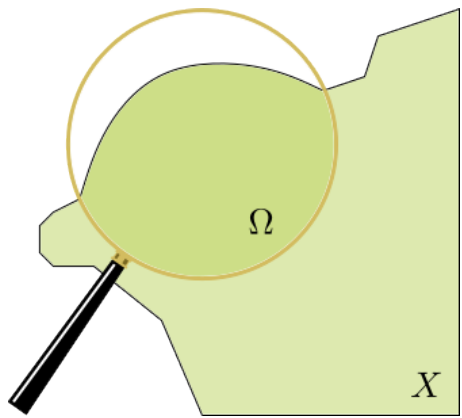

**Figure 1.2 - Definition of** $\text{Conv}_\kappa^o(\Omega)$ - Heuristically, we think of $\Omega$ as the aperture of a glass through which the convex set $X$ is examined, so that $\Omega$ is the region where the structure of $\partial X$ is better understood.

Finally, given a closed subset $Y$ of $\mathbb{H}^{m+1} \cup \partial_\infty \mathbb{H}^{m+1}$, let $\text{Hull}(Y)$ denote its convex hull in $\mathbb{H}^{m+1} \cup \partial_\infty \mathbb{H}^{m+1}$. The topological boundaries of convex hulls have precise structures which we now describe (see Section 4.5 of [18]). Note first that the convex hulls of finite subsets of $\mathbb{H}^{m+1} \cup \partial_\infty \mathbb{H}^{m+1}$ are precisely the, possibly ideal, polyhedra in this space. For this reason, we view convex hulls as generalized polyhedra, and just as the relative boundary of any polyhedron is a union of polyhedra of lower dimension, so too is the relative boundary of any convex hull the union of convex hulls of lower dimension. More precisely, given a closed subset $Y$ of $\mathbb{H}^{m+1} \cup \partial_\infty \mathbb{H}^{m+1}$, and a boundary point $x$ of $\text{Hull}(Y)$ lying away from $Y$, then every supporting totally geodesic hyperplane $H$ of $\text{Hull}(Y)$ at $x$ satisfies,

$$F_H(Y) := H \cap \text{Hull}(Y) = \text{Hull}(H \cap Y) . \tag{1}$$

We call $F_H(Y)$ the *face* of $\text{Hull}(Y)$ supported on $H$. Let $\text{Faces}(Y)$ denote the set of all faces of $\text{Hull}(Y)$.

We are now ready to state our main result.

**Theorem 1.1, Propagation of singularities**

*Let $K$ be a convex curvature function, let $\kappa > 0$ be a positive real number, let $\Omega$ be a convex, open subset of $\mathbb{H}^{m+1}$, and let $X$ be an element of $\partial \text{Conv}_\kappa^o(\Omega)$ with non-trivial interior. There exists a closed subset $Y$ of $\partial X \cap \partial \Omega$, and a subset $\mathcal{F}$ of $\text{Faces}(Y)$ such that*

$$\text{Sing}(X, \Omega) = \bigcup_{F \in \mathcal{F}} F \cap \Omega . \tag{2}$$

In other words, singularities of $\partial X$ in $\Omega$ propagate along entire faces towards the boundary of $\Omega$, where they can often then be eliminated by elementary geometric arguments, as we will see presently. Theorem 1.1 is proven in Section 5.

A nice application of this result concerns the study of the ideal Plateau problem for hypersurfaces of constant $K$-curvature in hyperbolic space. Indeed, we obtain a simpler, more geometric proof of the following result (c.f. [9] and [10]).

**Theorem 1.2, Ideal Plateau problem**

*Let $K$ be a convex curvature function, and let $\Omega \subseteq \partial_\infty \mathbb{H}^{m+1}$ be a proper open subset with $C^2$ boundary. For all $\kappa \in ]0, 1[$, there exists a convex subset $X \subseteq \mathbb{H}^{m+1} \cup \partial_\infty \mathbb{H}^{m+1}$ such that*

*(1) $X \cap \partial_\infty \mathbb{H}^{m+1} = \Omega^c$; and*

*(2) $\partial X \cap \mathbb{H}^{m+1}$ is smooth and of constant $K$-curvature equal to $\kappa$.*

**Remark 1.1.** Denoting $\Sigma := \partial X \cap \mathbb{H}^{m+1}$, we see that $\partial_\infty \Sigma = \partial \Omega$, so that the boundary of $X$ does indeed yield a solution of the ideal Plateau problem in the usual sense.

**Remark 1.2.** An elementary compactness argument shows that the principal curvatures of $\partial X$ are globally bounded above and below. Likewise, in the upper half-space model of $\mathbb{H}^{m+1}$, where $\partial_\infty \mathbb{H}^{m+1}$ identifies with $\mathbb{R}^m \cup \{\infty\}$, if $\Omega$ is taken to be a proper, open subset of $\mathbb{R}^m$, then an elementary convexity argument also shows that $\partial X \cap \mathbb{H}^{m+1}$ is a graph of some smooth function $u : \Omega \to [0, \infty[$.

Theorem 1.2 is proven in Section 6. This result in fact holds under a weaker hypothesis on the boundary of $\Omega$, not addressed in [9] and [10], known as the exterior ball condition. This will be discussed in Section 6.

Theorem 1.1 is a straightforward consequence of two ingredients. The first is an adaptation of the remarkable estimates obtained by Sheng–Urbas–Wang in [17], which build on the works [14] of Pogorelov and [12] and [13] of Ivochkina. This will be discussed in Section 4 where, for the reader's convenience, we will describe their construction in some detail. The second ingredient concerns the geometric properties of convex subsets of affine flat spaces, which also hold for convex subsets of hyperbolic space via the Kleinian parametrization (see, for example, [6]). These properties will be discussed in Section 5, where frequent reference will be made to our text [18]. The key result relating these two theories is Lemma 5.3, which transforms the analytic estimates obtained in Section 4 into local geometric properties of singular points.

**2 - Acknowledgements.** The author is grateful to Thierry Barbot for enlightening conversations, and to François Fillastre for helpful comments. This note was written whilst the author was benefitting from a CNRS Poste Rouge hosted by the Fédération de Recherche des Unités de Mathématiques de Marseille (FRUMAM) and the Laboratoire de Mathématiques of the Université d'Avignon.

**3 - Convex curvature functions.** We now introduce the class of curvature functions that will be of interest to us (c.f. [3]). Let $\Gamma^+ \subseteq \mathbb{R}^m$ denote the open cone of $m$-dimensional vectors all of whose components are strictly positive, and let $\overline{\Gamma}^+$ denote its closure in $\mathbb{R}^m$. We say that a function $K \in C^0(\overline{\Gamma}^+) \cap C^\infty(\Gamma^+)$ is a *convex curvature function* whenever it possesses the following properties.

(1) $K(x)$ is invariant under permutation of the components of $x$;

(2) $K$ is homogeneous of order $1$;

(3) $K(1, \cdots, 1) = 1$;

(4) $K(x) > 0$ for all $x \in \Gamma^+$ and $K(x) = 0$ for all $x \in \partial \Gamma^+$;

(5) for all $x \in \Gamma^+$, and for all $1 \leqslant i \leqslant m$,

$$(\partial_i K)(x) > 0; \text{ and} \qquad (3)$$

(6) $K$ is concave.

Convex curvature functions are used in the study of immersed hypersurfaces as follows. Let $Y := Y^{m+1}$ be an $(m+1)$-dimensional riemannian manifold, let $X := X^m$ be an $m$-dimensional manifold, let $e : X \to Y$ be an immersion, and let $\kappa_{1,e}, \cdots, \kappa_{m,e}$ denote its principal curvatures. We say that $e$ is *infinitesimally strictly convex* (ISC) whenever every one of its principal curvatures is strictly positive at every point. When this holds, we define its $K$-curvature $K_e : X \to \mathbb{R}$ by

$$K_e := K(\kappa_1, \cdots, \kappa_m) . \qquad (4)$$

The geometric significance of Properties 1 to 4 now becomes clear. Property 1 ensures that $K_e$ is well-defined, Property 2 ensures that it scales appropriately under rescalings of the ambient metric, Property 3 is a normalization condition according to which the unit sphere in euclidean space has unit curvature, and Property 4 ensures that the limit of a sequence of ISC hypersurfaces can only fail to be ISC at points where its $K$-curvature vanishes. In contrast, Properties 5 and 6 arise from analytic considerations. Indeed, Property 5 ensures that the Jacobi operator of $K$-curvature is elliptic, a fact that we will only use implicitly in what follows, whilst Property 6 is natural, and usually indispensable, with the current technology of non-linear analysis.

In what follows, it will be more convenient to study $K$-curvature in terms of the shape operator of the immersion, and thus to view $K$ as a function of symmetric matrics. Abusing terminology, we therefore denote also by $\Gamma^+$ the open cone of positive-definite matrices in $\mathrm{End}(\mathbb{R}^m)$, and by $K : \overline{\Gamma}^+ \to \mathbb{R}$ the function given by

$$K(A) := K(\lambda_1(A), \cdots, \lambda_m(A)) \ ,$$

where, for any symmetric matrix $A$, $\lambda_1(A), \cdots, \lambda_m(A)$ denote its eigenvalues. Note that Property 1 ensures that this new function is well-defined and invariant under the action of $\mathrm{O}(m)$ by conjugation on $\overline{\Gamma}^+$. An alternative approach consists of identifying $\mathbb{R}^m$ with the subspace of diagonal matrices in $\mathrm{End}(\mathbb{R}^m)$, and then viewing this new function as the unique $\mathrm{O}(m)$-invariant extension of the former. Any resulting ambiguity should be easily resolved by context.

Now let $A$ be a diagonal matrix with eigenvalues $\lambda_1 \geqslant \cdots \geqslant \lambda_m > 0$. $\mathrm{O}(m)$-invariance yields,

$$DK(A) \cdot M = \mathrm{Tr}(BM) \ ,$$

where $B$ is also diagonal, with respective eigenvalues

$$\mu_i := (\partial_i K)(\lambda_1, \cdots, \lambda_m) \ .$$

Furthermore, by Property 6 (concavity),

$$0 < \mu_1 \leqslant \cdots \leqslant \mu_m \ ,$$

and by Property 2 (homogeneity),

$$\sum_k \mu_k \lambda_k = DK(A) \cdot A = K(A) \ . \tag{5}$$

In what follows, it will also be useful to denote

$$\mu := \mathrm{Tr}(B) = \sum_k \mu_k \text{ and } \hat{\mu} := \mathrm{Tr}(BA^2) = \sum_k \mu_k \lambda_k^2 \ . \tag{6}$$

In particular, by (5) and Property 6 again,

$$\mu = K(A) + \sum_k \mu_k (1 - \lambda_k) \geqslant K(\mathrm{Id}) = 1 \ . \tag{7}$$

Finally (see, for example, Lemma 2.3 of [17] and Lemma 1.1 of [8]), the second derivative of $K$ at $A$ satisfies, for every symmetric matrix $M$,

$$DK(A)^{mn,pq} M_{mn} M_{pq} = \sum_{i,j} (\partial_i \partial_j K) M_{ii} M_{jj} + 2 \sum_{i>j} \frac{(\mu_i - \mu_j)}{(\lambda_i - \lambda_j)} M_{ij}^2 \ , \tag{8}$$

from which it follows, in particular, that this function is also concave.

**4 - Controlling the principal curvatures.** It will henceforth be convenient to view $\mathbb{H}^{m+1}$ as a quadric in $\mathbb{R}^{m+1,1}$, that is

$$\mathbb{H}^{m+1} := \left\{ x \in \mathbb{R}^{m+1,1} \mid \langle x, x \rangle_{m+1,1} = -1 \text{ and } x_{m+2} > 0 \right\} \ , \tag{9}$$

where $\langle \cdot, \cdot \rangle_{m+1,1}$ here denotes the lorentzian inner product with signature $(m+1, 1)$. In particular, every unit tangent vector of $\mathbb{H}^{m+1}$ is an element of $(m, 1)$-dimensional de Sitter space, that is

$$\mathrm{d}\mathbb{S}^{m,1} := \left\{ x \in \mathbb{R}^{m+1,1} \mid \langle x, x \rangle_{m+1,1} = 1 \right\} \ . \tag{10}$$

We will make considerable use of these identifications in what follows.

Let $K$ be a convex curvature function, let $X$ be an $m$-dimensional manifold, let $e : X \to \mathbb{H}^{m+1}$ be an ISC immersion of constant $K$-curvature equal to $\kappa > 0$, and let $\hat{e} : X \to \mathrm{d}\mathbb{S}^{m,1}$ denote its unit normal vector field, which we view as an immersion in its own right. Recall (see, for example, [7]) the natural duality of $\mathbb{H}^{m+1}$ and $d\mathbb{S}^{m,1}$ whereby the three fundamental forms of $\hat{e}$ are related to those of $e$ by

$$\mathrm{I}_{\hat{e}} = \mathrm{III}_e \ , \ \mathrm{II}_{\hat{e}} = \mathrm{II}_e \ , \ \mathrm{III}_{\hat{e}} = \mathrm{I}_e \ . \tag{11}$$

In particular, the principal curvatures of $\hat{e}$ are related to those of $e$ by

$$\kappa_{i,\hat{e}} = \frac{1}{\kappa_{i,e}} \ \forall 1 \leqslant i \leqslant m \ , \tag{12}$$

so that $\hat{e}$ is also an ISC immersion with constant $\hat{K}$-curvature equal to $1/\kappa$, where the curvature function $\hat{K}$ is given by

$$\hat{K}(x_1, \cdots, x_m) := K(1/x_1, \cdots, 1/x_m)^{-1} \ . \tag{13}$$

We will only make implicit use of this fact in what follows.

We aim to control the principal curvatures of $e$ in terms of the geometries of $e$ and $\hat{e}$, using a modification of the estimate [17] of Sheng–Urbas–Wang. In what follows, we aim to emphasize the symmetries of the roles played by the auxiliary functions $f$ and $\hat{f}$ that will be introduced presently.

Let $A := A_e$ denote the shape operator of $e$, and let $\lambda_{1,e} \geq \cdots \geq \lambda_{m,e} > 0$ denote its eigenvalues, that is, the principal curvatures of $e$. In what follows, we will use a semi-colon ";" to denote covariant differentiation with respect to the Levi-Civita covariant derivative of the metric induced by $e$. Note, in particular, that, by the Codazzi-Mainardi equations, for all $(i,j,k)$,

$$A_{ij;k} = A_{kj;i} \ , \tag{14}$$

so that $\nabla A$ is symmetric under every permutation of its components. For this reason, we will henceforth suppress the semi-colon in the first covariant derivative of $A$.

Now let $x_0$ be a point of $X$. Let $f$ be a positive, convex function defined in a neighbourhood $V$ of $e(x_0)$ in $\mathbb{H}^{m+1}$, and let $\hat{f}$ be a positive, concave function defined in a neighbourhood $\hat{V}$ of $\hat{e}(x_0)$ in $d\mathbb{S}^{m,1}$. For $\alpha > 0$ a positive constant to be determined presently, consider the function

$$\Phi_\alpha := \log(\lambda_1) + \alpha \log(f \circ e) - \log(\hat{f} \circ \hat{e}) \ . \tag{15}$$

We aim to estimate the $K$-laplacian of this function at $x_0$ whenever $x_0$ is a critical point.

Let $0 < \mu_{1,e} \leqslant \cdots \leqslant \mu_{m,e}$ denote the eigenvalues of $B_e := DK(A_e)$, as in Section 3. We define the $K$-laplacian of any smooth function $\phi : X \to \mathbb{R}$ by

$$\Delta^K \phi := \sum_k \mu_k \mathrm{Hess}^X(\phi)_{kk} \ , \tag{16}$$

where $\mathrm{Hess}^X$ denotes the Hessian operator of $X$ with respect to the metric induced by $e$. We require the following elementary formula for the $K$-laplacian of the logarithm of a positive function $\phi$.

$$\Delta^K \log(\phi) = \frac{1}{\phi} \Delta^K \phi - \frac{1}{\phi} \sum_{k=1}^m \mu_k \phi_k^2 \ . \tag{17}$$

We will also require the following *generalized Simons identity* (c.f. [15]).

**Lemma 4.1, Generalized Simons identity**

Let $Y$ be a riemannian space-form of constant curvature equal to $c$, let $e : X \to Y$ be a smooth codimension 1 immersion, and let $A$ denote its shape operator. Let $x$ be a point of $X$, let $e_1, \cdots, e_m \in T_x X$ be an orthonormal basis of eigenvectors of $A(x)$, and let $\lambda_1, \cdots, \lambda_m$ denote its corresponding eigenvalues. With respect to this basis, for all $(i, j)$,

$$A_{ii;jj} = A_{jj;ii} + (\lambda_i \lambda_j + c)(\lambda_i - \lambda_j) \ . \tag{18}$$

**Proof:** Indeed, let $\overline{R}$ denote the Riemann curvature tensor of $Y$ and let $R$ denote the Riemann curvature tensor of the metric that $e$ induces over $X$. Using the Codazzi-Mainardi equations (14), we find that, for all $(i, j)$,

$$A_{ii;jj} = A_{ji;ij} = A_{ji;ji} + R_{ijj}{}^p A_{pi} + R_{iji}{}^p A_{jp} = A_{jj;ii} + R_{ijj}{}^p A_{pi} + R_{iji}{}^p A_{jp} \ ,$$

where Einstein's summation convention is used in the index $p$, but not in the indices $i$ and $j$. By Gauss' equation,

$$R_{ijk}{}^l = \overline{R}_{ijk}{}^l + A_i^l A_{jk} - A_j^l A_{ik} \ ,$$

which, upon substitution into the preceding relation yields, for all $(i, j)$,

$$A_{ii;jj} = A_{jj;ii} + \overline{R}_{ijj}{}^p A_{pi} + \overline{R}_{iji}{}^p A_{jp} + A_i^p A_{jj} A_{pi} - A_j^p A_{ij} A_{pi} + A_i^p A_{ji} A_{jp} - A_j^p A_{ii} A_{jp}$$
$$= A_{jj;ii} + c\lambda_i - c\lambda_j + \lambda_i^2 \lambda_j - \delta_{ij} \lambda_i^3 + \delta_{ij} \lambda_i^3 - \lambda_j^2 \lambda_i$$
$$= A_{jj;ii} + (\lambda_i \lambda_j + c)(\lambda_i - \lambda_j) \ ,$$

as desired. $\square$

Finally, we recall (see [4]) that a continuous function $\phi$ satisfies the differential inequality

$$\Delta^K \phi \geqslant a$$

in the *weak sense* at some point $x$ whenever there exists a smooth function $\psi$ defined in a neighbourhood of this point such that

$$\phi \geqslant \psi \ , \ \phi(x) = \psi(x) \ , \text{ and } (\Delta^K \psi)(x) \geqslant a \ .$$

Note that this is stronger than satisfying the above differential inequality in the *viscosity sense* (see [5]), since the latter does not require the function $\psi$ to exist.

**Lemma 4.2**

When $e$ has constant K-curvature equal to $\kappa$, the function $\lambda_1$ satisfies

$$\Delta^K \log(\lambda_1) \geqslant \frac{2}{\lambda_1} \sum_{k>1} \frac{(\mu_k - \mu_1)}{(\lambda_1 - \lambda_k)} A_{11k}^2 - \frac{1}{\lambda_1^2} \sum_k \mu_k A_{11k}^2 - (\hat{\mu} + \mu) + \kappa \left( \lambda_1 + \frac{1}{\lambda_1} \right) \ , \tag{19}$$

*in the weak sense.*

**Remark 4.1.** Since we aim to apply a maximum principle, we will be interested in establishing mechanisms to ensure that the right-hand since of (19) is positive. We already note that, when $\mu_k \geqslant 2\mu_1$, the $k$-th summand of the second term on the right-hand side is dominated by the $k$-th summand of the first term. The case where $\mu_k < 2\mu_1$, as well as the third term on the right-hand side, will be addressed via the auxiliary functions $\log(f \circ e)$ and $\log(\hat{f} \circ \hat{e})$.

**Proof:** Let $x_0$ be a point of $X$, let $e_1, \cdots, e_m \in T_{x_0} X$ be an orthonormal basis of eigenvectors of $A(x_0)$, and denote also by $e_1, \cdots, e_m$ the orthonormal frame obtained by parallel transport of this basis along geodesics passing through $x_0$. Note that, with the frame defined in this manner, at $x_0$, for all $(i, j)$,

$$\nabla_{e_i} e_j = \nabla_{e_i} \nabla_{e_i} e_j = 0 \ . \tag{20}$$

Consider now the function
$$a_{11} := \langle Ae_1, e_1 \rangle \ .$$
Trivially $\lambda_1 \geqslant a_{11}$ near $x_0$, and $\lambda_1(x_0) = a_{11}(x_0)$, so that, by the above definition,
$$(\Delta^K \lambda_1)(x_0) \geqslant (\Delta^K a_{11})(x_0)$$
in the weak sense.

In order to estimate $(\Delta^K a_{11})(x_0)$, we consider, more generally, the matrix-valued function
$$a_{ij} := \langle Ae_i, e_j \rangle \ .$$
By (20), at $x_0$, for all $(i,j,k)$,
$$D_{e_k} a_{ij} = A_{ijk} \ , \quad D_{e_k} D_{e_k} a_{ij} = A_{ij;kk} \ , \text{ and } \text{Hess}^X(a_{11})(e_k, e_k) = D_{e_k} D_{e_k} a_{11} = A_{11;kk} \ . \tag{21}$$
Twice differentiating the identity $K(a) = \kappa$ yields, for all $k$,
$$\sum_{p,q} DK(a)^{pq} D_{e_k} D_{e_k} a_{pq} + \sum_{p,q,m,n} D^2 K(a)^{pq,mn} (D_{e_k} a_{pq})(D_{e_k} a_{mn}) = 0 \ ,$$
so that, at $x_0$, for all $k$,
$$\sum_{p,q} DK(a)^{pq} A_{pq;kk} = -\sum_{p,q,m,n} D^2 K(a)^{pq,mn} A_{pqk} A_{mnk} \ . \tag{22}$$
Using (5), (6), (21), (22), and the generalized Simons identity (18), we thus obtain
$$\Delta^K a_{11} = \sum_k \mu_k A_{11;kk}$$
$$= \sum_k \mu_k \big(A_{kk;11} + (\lambda_1 \lambda_k - 1)(\lambda_1 - \lambda_k)\big) \ ,$$
$$= -\sum_{p,q,m,n,k} D^2 K(\alpha)^{pq,mn} A_{pqk} A_{mnk} + \lambda_1^2 \kappa - \lambda_1 \mu - \lambda_1 \hat{\mu} + \kappa \ ,$$
so that, by (8),
$$\Delta^K a_{11} \geqslant 2 \sum_{k>1} \frac{(\mu_k - \mu_1)}{(\lambda_1 - \lambda_k)} A_{11k}^2 - \lambda_1(\hat{\mu} + \mu) + \kappa \lambda_1 \left(\lambda_1 + \frac{1}{\lambda_1}\right) \ ,$$
and the result now follows by (17). □

The calculations of the $K$-laplacians of the remaining two components of $\Phi_\alpha$ is more straightforward. We will use the following elementary but useful formula for hessians of restrictions of functions.

**Lemma 4.3**

Let $Y$ be a semi-riemannian manifold, let $X$ be a non-degenerate hypersurface in $Y$, let $\nu$ denote its unit normal vector field, and let $\epsilon$ denote the norm-squared of $\nu$. Given a point $x$ of $X$, and a $C^2$ function $\phi$ defined over some neighbourhood of this point in $Y$,
$$\text{Hess}^X(\phi) = \text{Hess}^Y(\phi)|_{TX} - \epsilon d\phi(\nu) \text{II}^X \ , \tag{23}$$
where $\text{II}^X$ here denotes the second fundamental form of $X$ with respect to the normal $\nu$.

**Proof:** Indeed, let $\nabla$ and $\overline{\nabla}$ denote respectively the Levi-Civita covariant derivatives of $X$ and $Y$. For all tangent vector fields $\xi$ and $\mu$ over $X$,
$$\text{Hess}^X(\phi)(\xi, \mu) = D_\mu D_\xi \phi - \phi(\nabla_\xi \mu)$$
$$= D_\mu D_\xi \phi - d\phi(\overline{\nabla}_\xi \mu - \epsilon \langle \overline{\nabla}_\xi \mu, \nu \rangle \nu)$$
$$= \text{Hess}^Y(\phi)(\xi, \mu) - \epsilon d\phi(\nu) \langle \mu, \overline{\nabla}_\xi \nu \rangle$$
$$= \text{Hess}^Y(\phi)(\xi, \mu) - \epsilon d\phi(\nu) \text{II}^X(\xi, \mu) \ ,$$
as desired. □

**Lemma 4.4**

Let $B > 0$ be such that
$$\text{Hess}^{\mathbb{H}^{m+1}}(f) \geqslant \frac{1}{B} f \text{Id} . \tag{24}$$

Then
$$\Delta^K \log(f \circ e) \geqslant \frac{1}{B} \mu - \frac{1}{f^2} \sum_k \mu_k f_k^2 - \frac{1}{f} df(\nu) \kappa . \tag{25}$$

**Proof:** Indeed, by (5), (6), (23) and (24),
$$\Delta^K (f \circ e) = \sum_k \mu_k \text{Hess}^X (f \circ e)_{kk} \geqslant \frac{1}{B} f \sum_k \mu_k - df(\nu) \sum_k \mu_k \lambda_k = \frac{1}{B} f\mu - df(\nu)\kappa ,$$

and the result now follows by (17). $\square$

**Lemma 4.5**

Let $B > 0$ be such that
$$\text{Hess}^{\text{d}\mathbb{S}^{m,1}}(\hat{f}) \leqslant -\frac{1}{B} \hat{f} \text{Id} . \tag{26}$$

Then
$$\Delta^K \log(\hat{f} \circ \hat{e}) \leqslant -\frac{1}{B} \hat{\mu} - \frac{1}{\hat{f}^2} \sum_k \mu_k \lambda_k^2 + \frac{1}{\hat{f}} \kappa d\hat{f}(\hat{\nu}) . \tag{27}$$

**Remark 4.2.** We understand (26) in the sense that
$$\text{Hess}^{\text{d}\mathbb{S}^{m,1}}(\hat{f})(\xi, \xi) \leqslant -\frac{1}{B} \hat{f} \langle \xi, \xi \rangle_{m+1,1}$$

for every spacelike vector $\xi$.

**Remark 4.3.** In fact, the $K$-laplacian of $e$ coincides with the $\hat{K}$-laplacian of $\hat{e}$, and the proof of Lemma 4.5 then reduces to an elementary calculation involving (23) almost identical to that used to prove Lemma 4.4. This is the main idea used below.

**Proof:** We will make use here of the natural duality between $\mathbb{H}^{m+1}$ and $\text{d}\mathbb{S}^{m,1}$ mentioned above. Since the shape operator $A_e$ of $e$ satisfies the Codazzi-Mainardi equations, we readily verify that the covariant derivative $\hat{\nabla}$ of $\text{I}_{\hat{e}} = \text{III}_e$ is related to that of $I_e$ by
$$\nabla_\xi \mu = \hat{\nabla}_\xi \mu + A^{-1}(\nabla_\xi A)\mu .$$

Consequently, for any function $\phi$ defined over $X$, and for all $(\xi, \mu)$,
$$\text{Hess}^X(\phi)(\xi, \mu) = \widehat{\text{Hess}}^X(\phi)(\xi, \mu) + d\phi(A^{-1}(\nabla_\xi A)\mu) ,$$

where here $\widehat{\text{Hess}}^X$ denotes the hessian operator of $\hat{\nabla}$. However, differentiating the identity $K(A) = \kappa$ yields, for all $i$,
$$\sum_k \mu_k A_{kki} = \text{Tr}(B \nabla_{e_i} A) = D_{e_i} K = 0 ,$$

so that, for any such $\phi$,
$$\Delta^K \phi = \sum_k \mu_k \text{Hess}^X(\phi)_{kk} = \sum_k \mu_k \widehat{\text{Hess}}^X(\phi)_{kk} .$$

Note now that the second fundamental form of $\hat{e}$ is equal to that of $e$. Since the unit normal $\hat{\nu}$ of $\hat{e}$ has negative norm-squared, using (5), (6), and (23), we therefore obtain
$$\Delta^K (\hat{f} \circ \hat{e}) = \sum_k \mu_k \widehat{\text{Hess}}^X(\hat{f} \circ \hat{e})_{kk} \leqslant -\frac{1}{B} \hat{f} \sum_k \mu_k \lambda_k^2 + d\hat{f}(\hat{\nu}) \sum_k \mu_k \lambda_k = -\frac{1}{B} \hat{f} \hat{\mu} + d\hat{f}(\hat{\nu}) \kappa ,$$

and the result now follows by (17). $\square$

We are now ready to provide our estimate for the $K$-laplacian of $\Phi_\alpha$.

**Lemma 4.6**

Let $K \subseteq \mathbb{H}^{m+1}$ and $\hat{K} \subseteq d\mathbb{S}^{m,1}$ be respectively compact subsets of $\mathbb{H}^{m+1}$ and $d\mathbb{S}^{m,1}$. Suppose that $e(X) \subseteq K$ and $\hat{e}(X) \subseteq \hat{K}$, and that there exists $B > 0$ such that, over $K$ and $\hat{K}$,

$$\text{Hess}(f) \geqslant \frac{1}{B} f \,, \hat{f} \geqslant \frac{1}{B} \,, \text{ and Hess}(\hat{f}) \leqslant -\left(1 + \frac{1}{B}\right)\hat{f} \,. \tag{28}$$

For all $\alpha > \text{Max}(2, B)$, there exists $C > 0$, which only depends on $\alpha$, $K$, $\hat{K}$, and $B$, with the following property. At any critical point $x_0$ of $\Phi_\alpha$ with $\Phi_\alpha(x_0) > C$,

$$\Delta^K \Phi_\alpha > 0 \tag{29}$$

in the weak sense.

**Proof:** Note first that $\nu_e := \hat{e}$ and $\nu_{\hat{e}} := e$ are respectively the unit normal vector fields over $e$ and $\hat{e}$. Since the images of these functions are relatively compact, upon increasing $B$ is necessary, we may suppose that

$$|df(\nu_e)|, |d\hat{f}(\nu_{\hat{e}})| < B \,,$$

and that, for all $k$,

$$|f_k|, |\hat{f}_k| < B \,.$$

It now follows by (19), (25), and (27) that

$$\Delta^K \Phi_\alpha \geqslant \kappa \lambda_1 + \frac{2}{\lambda_1} \sum_{k>1} \frac{(\mu_k - \mu_1)}{(\lambda_1 - \lambda_k)} A_{11k}^2 - \frac{1}{\lambda_1^2} \sum_k \mu_k A_{11k}^2$$
$$- \mu + \frac{\alpha}{B}\mu - \frac{\alpha}{f^2} \sum_k \mu_k f_k^2$$
$$- \hat{\mu} + \left(1 + \frac{1}{B}\right)\hat{\mu} + \frac{1}{\hat{f}^2} \sum_k \mu_k \lambda_k^2 \hat{f}_k^2 - O\left(1 + \frac{\alpha}{f}\right) \,,$$

so that, for $\alpha \geqslant B$,

$$\Delta^K \Phi \geqslant \kappa \lambda_1 - O\left(1 + \frac{\alpha}{f}\right) + \sum_k \left(P_{1,k} + P_{2,k} - N_{1,k} - N_{2,k}\right) \,,$$

where the summands

$$P_{1,k} := \frac{2}{\lambda_1} \frac{(\mu_k - \mu_1)}{(\lambda_1 - \lambda_k)} A_{11k}^2 \,, \text{ and}$$

$$P_{2,k} := \frac{1}{\hat{f}^2}\left(1 + \frac{1}{B^5}\right)\mu_k \lambda_k^2 \hat{f}_k^2$$

yield non-negative contributions, whilst the summands

$$N_{1,k} := \frac{1}{\lambda_1^2} \mu_k A_{11k}^2 \,, \text{ and}$$

$$N_{2,k} := \frac{\alpha}{f^2} \mu_k f_k^2$$

yield non-positive contributions (by convention, we have defined $P_{1,1} := 0$).

We will show that, for all $k$, the terms $N_{1,k}$ and $N_{2,k}$ are dominated by the terms $P_{1,k}$ and $P_{2,k}$ provided that $\alpha$ and $\Phi_\alpha$ are sufficiently large. To this end, we observe first that, at any critical point $x_0$ of $\Phi_\alpha$, for all $k$,

$$\frac{1}{\lambda_1} A_{11k} + \frac{\alpha f_k}{f} - \frac{\lambda_k \hat{f}_k}{\hat{f}} = 0 \,. \tag{30}$$

We observe also that, since all principal curvatures are positive, by (5),

$$\mu_1 \lambda_1 \leqslant \sum_k \mu_k \lambda_k = \kappa ,$$

so that

$$\mu_1 \leqslant \frac{\kappa}{\lambda_1} . \tag{31}$$

Now choose $0 < \epsilon < 1$. The desired dominations of $N_{1,k}$ and $N_{2,k}$ will be obtained upon using (30) in two different ways, depending on the size of $\mu_k$. In the first case, when $\mu_k < 2\kappa/(1-\epsilon)\lambda_1$, we use (30) to obtain

$$\frac{1}{\lambda_1^2} A_{11k}^2 \leqslant \frac{(1+B^5)\alpha^2 f_k^2}{f^2} + \frac{1}{\hat{f}^2}\left(1 + \frac{1}{B^5}\right)\lambda_k^2 \hat{f}_k^2 ,$$

which yields

$$P_{2,k} - N_{1,k} = \frac{1}{\hat{f}^2}\left(1 + \frac{1}{B^5}\right)\mu_k \lambda_k^2 \hat{f}_k^2 - \frac{1}{\lambda_1^2}\mu_k A_{11k}^2 \geqslant -\frac{\mu_k(1+B^5)\alpha^2 f_k^2}{f^2} = -\mathrm{O}\left(\frac{\alpha^2}{\lambda_1 f^2}\right) ,$$

so that, for $\alpha > 1$,

$$P_{1,k} + P_{2,k} - N_{1,k} - N_{2,k} \geqslant -\mathrm{O}\left(\frac{\alpha^2}{\lambda_1 f^2}\right) .$$

In the second case, when $\mu_k \geqslant 2\kappa/(1-\epsilon)\lambda_1$, we use (30) to obtain

$$\frac{\alpha f_k^2}{f^2} \leqslant \frac{2}{\alpha \lambda_1^2} A_{11k}^2 + \frac{2}{\alpha \hat{f}^2}\lambda_k^2 \hat{f}_k^2 ,$$

which yields, for $\alpha > 2$,

$$P_{2,k} - N_{1,k} - N_{2,k} = \frac{1}{\hat{f}^2}\left(1 + \frac{1}{B^5}\right)\mu_k \lambda_k^2 \hat{f}_k^2 - \frac{1}{\lambda_1^2}\mu_k A_{11k}^2 - \frac{\alpha \mu_k f_k^2}{f^2}$$

$$\geqslant -\frac{1}{\lambda_1^2}\mu_k\left(1 + \frac{2}{\alpha}\right)A_{11k}^2 .$$

Bearing in mind (31), this yields

$$P_{2,k} - N_{1,k} - N_{2,k} \geqslant -\frac{2(\mu_k - \mu_1)}{\lambda_1(\lambda_1 - \lambda_k)}\frac{(1 + 2/\alpha)}{(1+\epsilon)}A_{11k}^2 ,$$

so that

$$P_{1,k} + P_{2,k} - N_{1,k} - N_{2,k} \geqslant \frac{2(\mu_k - \mu_1)}{\lambda_1(\lambda_1 - \lambda_k)}\frac{(\epsilon - 2/\alpha)}{(1+\epsilon)}A_{11k}^2 .$$

Combining these relations, we find that, for $\epsilon > 2/\alpha$, there exist $B'' \geqslant B' \geqslant B > 0$ such that

$$\Delta^K \Phi \geqslant \lambda_1\left(\kappa - \frac{B'}{\lambda_1 f^2}\right) \geqslant \lambda_1\left(\kappa - \frac{B''}{\lambda_1 f^\alpha}\right) \geqslant \lambda_1(\kappa - (B'')^2 e^{-\Phi}) ,$$

and the result now follows upon setting $C := \log((B'')^2/\kappa)$. $\square$

On the propagation of singularities...

**5 - Propagation of singularities.** We now use Lemma 4.6 to show how the singularities of ISC hypersurfaces of constant $K$-curvature propagate along totally geodesic hypersurfaces. We limit ourselves to hypersurfaces which bound convex sets, and we recall the formalism given in the introduction. Let Conv denote the space of closed, convex subsets of $\mathbb{H}^{m+1}\cup\partial_\infty\mathbb{H}^{m+1}$, furnished with the Hausdorff topology. Recall (see, for example, [6]), that the Kleinian parametrization maps $\mathbb{H}^{m+1}\cup\partial_\infty\mathbb{H}^{m+1}$ to the closed unit ball $\overline{\mathbb{B}}^{m+1}$ in $\mathbb{R}^{m+1}$, sending geodesics to straight lines, so that Conv identifies with the space of closed, convex subsets of $\overline{\mathbb{B}}^{m+1}$. In particular, Conv is compact (see, for example, Lemma 4.3 of [18]). Let $\text{Conv}^o$ denote the subspace of Conv consisting of those closed, convex subsets which have non-trivial interior. Given an open subset $\Omega$ and a positive real number $\kappa$, let $\text{Conv}^o_\kappa(\Omega)$ denote the subspace consisting of those convex subsets $X$, with non-trivial interior, such that $\partial X\cap\Omega$ is smooth and of constant $K$-curvature equal to $\kappa$ (see Figure 1.2).

We will be interested in the topological boundary of $\text{Conv}^o_k(\Omega)$ in $\text{Conv}^o$. Given an element $X\in\partial\text{Conv}^o_k(\Omega)$ with non-trivial interior, and a point $x\in\partial X\cap\Omega$, we will say that $X$ is *regular* at $x$ whenever $X\in\text{Conv}^o_k(\Omega')$, for some smaller neighbourhood $\Omega'$ of $x$, and we will say that it is *singular* at this point otherwise. Let $\text{Reg}(X,\Omega)$ denote the set of points of $\partial X\cap\Omega$ at which $X$ is regular, and let $\text{Sing}(X,\Omega)$ denote the set of points at which it is singular. Note that, by definition, $\text{Reg}(X,\Omega)$ is an open subset of $\partial X\cap\Omega$, whilst $\text{Sing}(X,\Omega)$ is closed. Finally, given a closed subset $Y$ of $\mathbb{H}^{m+1}\cup\partial_\infty\mathbb{H}^{m+1}$, let $\text{Hull}(Y)$ denote its convex hull in $\mathbb{H}^{m+1}\cup\partial_\infty\mathbb{H}^{m+1}$, and let $\text{Faces}(Y)$ denote the set of all its faces.

The remainder of this section will be devoted to proving Theorem 1.1.

**Theorem 1.1**

*Let $K$ be a convex curvature function, let $\kappa>0$ be a positive real number, let $\Omega$ be a convex, open subset of $\mathbb{H}^{m+1}$, and let $X$ be an element of $\partial\text{Conv}^o_\kappa(\Omega)$ with non-trivial interior. There exists a closed subset $Y$ of $\partial X\cap\partial\Omega$, and a subset $\mathcal{F}$ of $\text{Faces}(Y)$ such that*

$$\text{Sing}(X,\Omega)=\bigcup_{F\in\mathcal{F}} F\cap\Omega\ . \tag{32}$$

This theorem will follow as a straightforward consequence of Lemma 4.6 together with the following characterization of convex hulls.

**Definition 5.1**

*Let $X\in\text{Conv}^o$ be a convex subset of $\mathbb{H}^{m+1}$ with non-trivial interior. We say that $X$ satisfies the* local geodesic property *(LGP) at a point $x\in\partial X$ whenever there exists an open geodesic segment $\Gamma$ such that*

$$x\in\Gamma\subseteq\partial X\ . \tag{33}$$

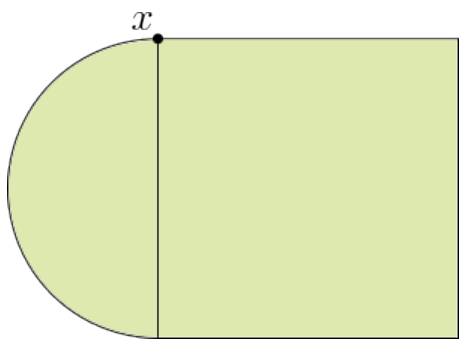

**Figure 5.3 - The local geodesic property** - Although the local geodesic property is not satisfied at the point $x$, this set is non-strictly convex at this point.

**Remark 5.1.** The set $X$ is trivially non-strictly convex at any point where it satisfies the LGP, but the converse is not true. Indeed, consider the union of a square and a semi-disk, as in Figure 5.3. This set is non-strictly convex at $x$ since its unique supporting tangent at this point also contains one side of the square, but it trivially does not satisfy the LGP at this point.

The LGP characterizes convex hulls, as the following result shows (see Theorem 4.8 of [18]).

**Lemma 5.2**

*Let $X$ be a compact, convex subset of $\mathbb{R}^{m+1}$ with non-trivial interior. Let $Y \subseteq \partial X$ denote the set of points of $\partial X$ at which $X$ satisfies the LGP, and let $\partial Y$ denote its topological boundary. Then*

$$Y = \partial X \cap \mathrm{Hull}(\partial Y) \ .$$

**Remark 5.2.** Since the Kleinian parametrization sends $\mathbb{H}^{m+1} \cup \partial_\infty \mathbb{H}^{m+1}$ to the closed unit ball $\overline{\mathbb{B}}^{m+1}$ in $\mathbb{R}^{m+1}$, sending geodesic segments to geodesic segments, this result also applies to the closed, convex subsets studied here.

By Lemma 5.2, Theorem 1.1 will follow as soon as we show that the singular points of $X$ in $\Omega$ are precisely those points at which it satisfies the LGP.

**Lemma 5.3**

*Let $\Omega$ by an open subset of $\mathbb{H}^{m+1}$, let $X$ be an element of $\mathrm{Conv}_k^o(\Omega)$ with non-trivial interior, and let $x_0$ be a point of $\partial X \cap \Omega$. $X$ is singular at $x_0$ if and only if it satisfies the LGP at this point.*

**Proof:** Note first that if $X$ satisfies the LGP at $x_0$, then it is singular at this point, for otherwise the lowest principal curvature of $\partial X \cap \Omega$ would vanish at this point, and therefore so too would the $K$-curvature of this hypersurface, which is absurd. Suppose now that $X$ does not satisfy the LGP at $x_0$. We will show that $X$ is regular at this point.

First note that, since $X$ has non-trivial interior, its set of supporting normals at $x_0$ is contained in an open hemisphere of $T^1_{x_0}\mathbb{H}$ and that there therefore exists a supporting normal $\nu_0 \in \mathrm{d}\mathbb{S}^{m,1}$ at $x_0$ and $\epsilon > 0$ such that, for every other supporting normal $\nu$ of $X$ at this point,

$$\langle \nu_0, \nu \rangle > 4\epsilon \ .$$

By continuity (see Lemma 4.4 of [18]), there exists a neighbourhood $U$ of $x_0$ in $\mathbb{H}^{m+1}$ such that, for all $x \in \partial X \cap V$, and for every supporting normal $\nu \in \mathrm{d}\mathbb{S}^{m,1}$ of $X$ at $x$,

$$\langle \nu_0, \nu \rangle > 3\epsilon \ .$$

We define $\hat{f} : \mathrm{d}\mathbb{S}^{m,1} \to \mathbb{R}$ and $\hat{V} \subseteq \mathrm{d}\mathbb{S}^{m,1}$ by

$$\hat{f}(y) = \langle \nu_0, y \rangle - \epsilon \text{ and}$$
$$\hat{V} = \hat{f}^{-1}(]2\epsilon, +\infty[) \ ,$$

and we claim that,

$$\mathrm{Hess}^{\mathrm{d}\mathbb{S}^{m,1}}(\hat{f}) = -(\hat{f} + \epsilon)\mathrm{Id} \ .$$

Indeed, let $N$ denote the outward-pointing unit normal vector field over $\mathrm{d}\mathbb{S}^{m,1}$ in $\mathbb{R}^{m+1,1}$. Since $N$ is spacelike, it follows by (23) that

$$\mathrm{Hess}^{\mathrm{d}\mathbb{S}^{m,1}}(\hat{f}) = \mathrm{Hess}^{\mathbb{R}^{m+1,1}}(\hat{f})|_{\mathrm{d}\mathbb{S}^{m,1}} - d\hat{f}(N)\mathrm{II}^{\mathrm{d}\mathbb{S}^{m,1}} = -(\hat{f} + \epsilon)\mathrm{Id} \ ,$$

as asserted.

Now let $V_0 \subseteq \mathbb{H}^{m+1}$ denote the unique open half-space such that $\partial V_0$ passes through $x_0$ with inward-pointing normal $\nu_0$ at this point. For all $r > 0$, let $\overline{B}_r(x_0)$ denote the closed ball of radius $r$ about $x_0$ in $\mathbb{H}^{m+1}$. Let $r > 0$ be such that $\overline{B}_r(x_0) \subseteq U$. Since $X$ does not satisfy the LGP at $x_0$, $\partial X \cap \partial V_0 \cap \partial B_r(x_0)$ is contained in an open hemisphere of $\partial B_r(x_0)$, bounded by a totally geodesic hypersurface $H_0 \subseteq \partial V_0$, say (see Theorem 4.19 and Lemma 4.21 of [18]). Rotating $V_0$ slightly about $H_0$ and then translating slightly away from $\nu_0$ thus yields an open half-space $V$, containing $x_0$, such that

$$\partial X \cap \overline{V} \subseteq U \ ,$$

On the propagation of singularities...

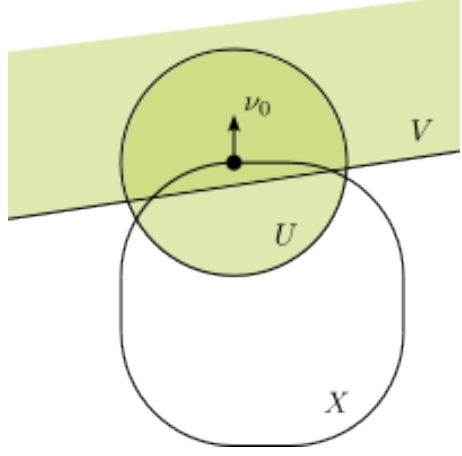

**Figure 5.4 - Constructing a suitable half-space** - Upon rotating slightly, and then translating slightly a suitable supporting open half-space of $X$ at $x_0$, we obtain an open half-space $V$, containing $x_0$ whose intersection with $X$ is a relatively compact subset of $U$.

and that this intersection is compact (see Figure 5.4). Recall now that $\partial V$ is the intersection of $\mathbb{H}^{m+1}$ with a non-degenerate, linear hyperplane $P \subseteq \mathbb{R}^{m+1,1}$ of mixed signature. Let $\nu_1 \in \mathbb{R}^{m+1,1}$ denote the unit vector normal to this hyperplane which points into $V$, and define $f : \mathbb{H}^{m+1} \to \mathbb{R}$ by

$$f(x) := \langle x, \nu_1 \rangle .$$

Since the unit normal vector over $\mathbb{H}^{m+1}$ in $\mathbb{R}^{m+1,1}$ is timelike, using (23), we obtain, as before,

$$\operatorname{Hess}^{\mathbb{H}^{m+1}}(f) = f\operatorname{Id} .$$

Now let $(X_m)_{m\in\mathbb{N}}$ be a sequence of elements of $\operatorname{Conv}^o_\kappa(\Omega)$ converging in the Hausdorff sense to $X$. We claim that $(\partial X_m \cap V)_{m\in\mathbb{N}}$ converges in the $C^\infty_{\text{loc}}$ sense towards $\partial X \cap V$, from which the regularity of $X$ at $x_0$ will immediately follow. Indeed, for all $m$, let $\nu_m$ denote the outward-pointing, unit, normal vector field over $\partial X_m \cap V$, let $A_m$ denote its shape operator, and, for all $\alpha > 0$, denote

$$\Phi_{\alpha,m} := \log(\|A_m\|) + \alpha \log(f) - \log(\hat{f} \circ \nu_m) .$$

Since $(X_m)_{m\in\mathbb{N}}$ converges in the Hausdorff sense to $X$, we may suppose that, for all $m$,

$$\partial X_m \cap \overline{V} \subseteq U ,$$

and that this intersection is compact. Furthermore, by Lemma 4.4 of [18], we may also suppose that, for all $m$,

$$\nu_m(\partial X_m \cap \overline{V}) \subseteq \hat{V} .$$

Note now that, for all $m$, $\Phi_{\alpha,m}$ converges to $-\infty$ near the boundary of $\partial X_m \cap \overline{V}$, and therefore attains its maximum at some interior point of this set. With $\alpha$ as in Lemma 4.6, it then follows by the maximum principle that there exists $C > 0$ such that, for all $m$,

$$\|A_m\| f^\alpha (\hat{f} \circ \nu_m)^{-1} = \operatorname{Exp}(\Phi_{\alpha,m}) \leqslant C ,$$

so that $(\|A_m\|)_{m\in\mathbb{N}}$ is locally uniformly bounded over $\partial X_m \cap V$. It now follows by elliptic regularity that $(\partial X_m \cap V)_{m\in\mathbb{N}}$ converges to $\partial X \cap V$ in the $C^\infty_{\text{loc}}$ sense, as asserted. This completes the proof. $\square$

We now complete the proof of Theorem 1.1.

**Proof of Theorem 1.1:** Indeed, by Lemma 5.3, $\mathrm{Sing}(X,\Omega)$ is precisely the set of points of $\partial X \cap \Omega$ at which $X$ satisfies the LGP. It follows by Lemma 5.2 that

$$\mathrm{Sing}(X,\Omega) = \partial X \cap \mathrm{Hull}(Y) \cap \Omega ,$$

where

$$Y := \partial \mathrm{Sing}(X,\Omega) \subseteq \partial \Omega .$$

Now let $x$ be a point of $\mathrm{Sing}(X,\Omega)$, and let $H$ be a supporting totally geodesic hypersurface of $X$ at $x$. Since $H$ is also a supporting totally geodesic hypersurface of $\mathrm{Hull}(Y)$ at $x$,

$$x \in F_H(Y) \cap \Omega = H \cap \mathrm{Hull}(Y) \cap \Omega \subseteq \partial X \cap \Omega .$$

It follows that $\mathrm{Sing}(X,\Omega)$ is indeed the union of intersections with $\Omega$ of faces of $\mathrm{Hull}(Y)$, and this completes the proof. $\square$

**6 - The asymptotic Plateau problem.** We conclude this note by applying Theorem 1.1 to the study of the asymptotic Plateau problem for convex curvature functions in $\mathbb{H}^{m+1}$. We require the following definition.

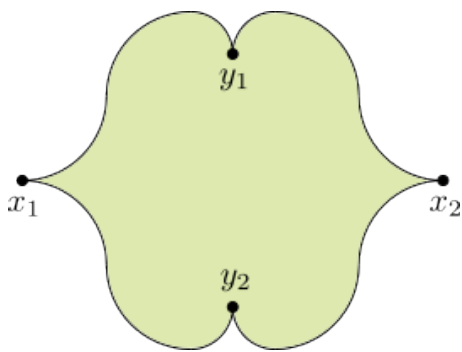

**Figure 6.5 - The exterior ball condition** - This set satisfies the exterior ball condition at every point of its boundary where it is smooth. It also satisfies the exterior ball condition at the non-smooth points $x_1$ and $x_2$, but not $y_1$ and $y_2$.

**Definition 6.1**

Let $\Omega$ be a proper open subset of $\partial_\infty \mathbb{H}^{m+1}$, and let $x$ be a boundary point of this set. We say that $\Omega$ satisfies the *exterior ball condition (EBC)* at this point whenever there exists a closed round ball $B \subseteq \partial_\infty \mathbb{H}^{m+1}$, lying in the complement of $\Omega$, which contains $x$ in its boundary (see Figure 6.5).

We now obtain a simpler proof of the following result of Guan–Spruck and Guan–Spruck–Szapiel (see [9] and [10]).

**Theorem 6.2**

Let $K$ be a convex curvature function, and let $\Omega \subseteq \partial_\infty \mathbb{H}^{m+1}$ be a proper open subset satisfying the EBC at every point of its boundary. For all $\kappa \in ]0,1[$, there exists a convex subset $X \subseteq \mathbb{H}^{m+1} \cup \partial_\infty \mathbb{H}^{m+1}$ such that

(1) $X \cap \partial_\infty \mathbb{H}^{m+1} = \Omega^c$; and

(2) $\partial X \cap \mathbb{H}^{m+1}$ is smooth and of constant $K$-curvature equal to $\kappa$.

**Remark 6.1.** Theorem 6.2 is actually slightly stronger than the results of [9] and [10], since the latter only address the case where the boundary is $C^2$.

**Proof:** We use the upper half-space model of $\mathbb{H}^{m+1}$, namely

$$\mathbb{H}^{m+1} := \{(x_1, \cdots, x_{m+1}) \mid x_{m+1} > 0\} ,$$

and we identify $\partial_\infty \mathbb{H}^{m+1}$ with $\mathbb{R}^m \cup \{\infty\}$. In particular, since $\Omega$ is a proper open subset of $\partial_\infty \mathbb{H}^{m+1}$ which satisfies the EBC at every point of its boundary, we may view it as a relatively compact open subset of $\mathbb{R}^m$. Let $(\Omega_i)_{i \in \mathbb{N}}$ be a sequence of open subsets of $\mathbb{R}^m$ with smooth boundary such that,

$$\Omega := \bigcup_{i \in \mathbb{N}} \Omega_i \ ,$$

and such that, for all $i$,

$$\overline{\Omega}_i \subseteq \Omega_{i+1} \ .$$

For all $i \in \mathbb{N}$, let $B_i$ denote the open horoball

$$B_i := \{(x_1, \cdots, x_{m+1}) \mid x_{m+1} > 1/i\} \ ,$$

and let $\hat{\Omega}_i$ denote the set of points of $\partial B_i$ lying directly above $\Omega_i$. For all $i$, solving the finite Plateau problem (see, for example, Theorem 1.2 of [19] or Theorem 1.3 of [10]) yields a closed, convex subset $X_i \in \mathrm{Conv}^o_\kappa(B_i)$ such that

$$X_i \subseteq \overline{B}_i \text{ and } X_i \cap \partial B_i = \hat{\Omega}_i^c \ .$$

By compactness of Conv, there then exists a closed, convex subset $X \in \mathrm{Conv}$ towards which $(X_i)_{i \in \mathbb{N}}$ subconverges in the Hausdorff sense. It will suffice to show that $X$ has the desired properties.

Before proceeding, it is worth recalling the structure of totally umbilic subspaces of $\mathbb{H}^{m+1}$ in the upper half-space model (c.f. [1]). In particular, this will yield a clear picture of how smoothness almost immediately follows from the EBC. In the upper half-space model, the intersections of hyperplanes and hyperspheres of $\mathbb{R}^{m+1}$ with $\mathbb{H}^{m+1}$ yield totally umbilic hypersurfaces of $\mathbb{H}^{m+1}$, and, furthermore, all totally umbilic hypersurfaces of $\mathbb{H}^{m+1}$ are accounted for in this manner. In addition, if $H$ is a hyperplane or a hypersphere which intersects the boundary $\mathbb{R}^m$ non-trivially, making an angle of $\theta \in [0, \pi/2[$ with the vertical axis along this intersection, then the shape operator of $H \cap \mathbb{H}^{m+1}$ is equal to $\sin(\theta)\mathrm{Id}$, so that its $K$-curvature is equal to $\sin(\theta)$. Finally, for any such $H$, $H \cap \mathbb{H}^{m+1}$ bounds a convex subset of $\mathbb{H}^{m+1}$. In the case where $H$ is a hyperplane, its outward-pointing normal points towards $\mathbb{R}^m$, in the case where it is a hypersphere with centre lying below $\mathbb{R}^m$, its outward-pointing normal points downwards, and in the case where it is a hypersphere with centre lying above $\mathbb{R}^m$, its outward-pointing normal points upwards.

We now show that $X \cap \partial_\infty \mathbb{H}^{m+1} = \Omega^c$. Indeed, trivially

$$\Omega^c \subseteq X \cap \partial_\infty \mathbb{H}^{m+1} \ .$$

Conversely, choose $x \in \Omega$, and let $B$ be an open ball in $\mathbb{R}^m$ such that

$$x \in B \subseteq \overline{B} \subseteq \Omega \ ,$$

and note that, for sufficiently large $i$,

$$\overline{B} \subseteq \Omega_i \ .$$

Now let $\hat{B}$ denote the open ball in $\mathbb{R}^{m+1}$ with centre lying below $\mathbb{R}^m$, whose intersection with $\mathbb{R}^m$ is equal to $B$, and which makes an angle of $\arcsin(\kappa)$ with the vertical axis along $\partial B$. By the geometric maximum principle, for all $i$,

$$X_i \cap \hat{B} = \emptyset \ .$$

Upon taking limits, it follows that

$$X \cap \hat{B} = \emptyset \ ,$$

so that $x \notin X$. Since $x \in \Omega$ is arbitrary, it follows that the intersection of $X \cap \partial_\infty \mathbb{H}^{m+1}$ with $\Omega$ is trivial, so that $X \cap \partial_\infty \mathbb{H}^{m+1} = \Omega^c$, as asserted.

We now show that $X \in \mathrm{Conv}^0_\kappa(\mathbb{H}^{m+1})$. To this end, we first show that $\partial X \cap \mathbb{H}^{m+1}$ makes an angle of $\theta := \arcsin(\kappa)$ with the vertical axis along $\partial \Omega$. Indeed, let $y$ be a point of $\partial \Omega$. By the EBC, there exists an open ball $B$ in $\mathbb{R}^m$ contained in the complement of $\Omega$ which contains $y$ in its closure. For all $r \in ]0, 1]$, let $B_r$ denote the open ball in $\mathbb{R}^m$ obtained upon contracting $B$ by a factor $r$ about its centre, and let $\hat{B}_r$ denote the open ball in $\mathbb{R}^{m+1}$ with centre lying above $\mathbb{R}^m$, whose intersection with $\mathbb{R}^m$ is equal to $B_r$, and which makes an angle of $\theta$ with the vertical axis along $\partial B_r$. By the geometric maximum principle again, for all $r$, and for all sufficiently large $i$,
$$\hat{B}_r \cap B_i \subseteq X_i \ .$$
Upon taking limits, first in $i$, and then in $r$, it follows that
$$\hat{B} \cap \mathbb{H}^{m+1} \subseteq X \ ,$$
so that, indeed, $\partial X \cap \mathbb{H}^{m+1}$ makes an angle of $\theta$ with $\mathbb{R}^m$ at $y$.

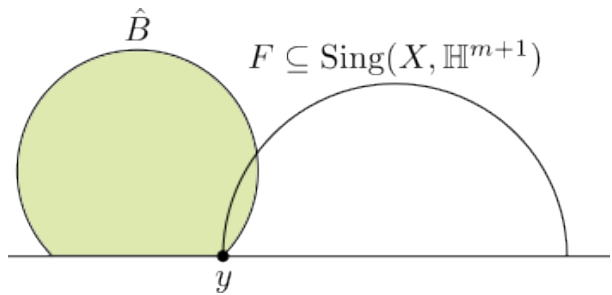

**Figure 6.6 - Towards a contradiction** - The set $\mathrm{Sing}(X,\mathbb{H}^{m+1})$ meets $\mathbb{H}^{m+1}$ orthogonally at $y$. However, by the geometric maximum principle, its intersection with the region marked in green is trivial, which is absurd.

It remains only to show that $\mathrm{Sing}(X,\mathbb{H}^{m+1})$ is empty. However, by Theorem 1.1, there exists a closed subset $Y$ of $\partial\Omega$ and a subset $\mathcal{F} \subseteq \mathrm{Faces}(Y)$ such that
$$\mathrm{Sing}(X,\mathbb{H}^{m+1}) = \bigcup_{F \in \mathcal{F}} F \cap \mathbb{H}^{m+1} \ .$$

Choose $x \in \mathrm{Sing}(X,\mathbb{H}^{m+1})$, let $F \in \mathcal{F}$ be a face of $\mathrm{Hull}(Y)$ containing this point, let $H$ denote the totally geodesic hypersurface on which it lies, and let $y \in (\partial_\infty H) \cap Y$ an extremity of $F$. Note, in particular, that $H$, and therefore also $F$, meets $\mathbb{R}^m$ orthogonally at $y$. However, with $\hat{B}$ as above, by the preceeding discussion,
$$F \cap \mathbb{H}^{m+1} \cap \hat{B} \subseteq \partial X \cap \mathbb{H}^{m+1} \cap \hat{B} = \emptyset \ ,$$
which is absurd (see Figure 6.6). It follows that $\mathrm{Sing}(X,\mathbb{H}^{m+1})$ is indeed trivial, so that $X \in \mathrm{Conv}^0_\kappa(\mathbb{H}^{m+1})$, and this completes the proof. $\square$

## 7 - Bibliography.